\documentclass[a4paper,11pt]{article}
\title{\bf{Fourier-Mukai transforms and canonical divisors}
}
\date{}
\author{Yukinobu Toda}

\usepackage{makeidx}

\usepackage{latexsym}
\usepackage{amscd}
\usepackage{amsmath}
\usepackage{amssymb}

\usepackage{amsthm}
\usepackage{float}
\usepackage[all]{xy}
\usepackage{graphicx}

\newcommand{\aA}{\mathcal{A}}
\newcommand{\bB}{\mathcal{B}}

\newcommand{\eE}{\mathcal{E}}
\newcommand{\fF}{\mathcal{F}}
\newcommand{\gG}{\mathcal{G}}
\newcommand{\hH}{\mathcal{H}}

\newcommand{\lL}{\mathcal{L}}

\newcommand{\oO}{\mathcal{O}}
\newcommand{\pP}{\mathcal{P}}
\newcommand{\qQ}{\mathcal{Q}}
\newcommand{\rR}{\mathcal{R}}

\newcommand{\tT}{\mathcal{T}}
\newcommand{\uU}{\mathcal{U}}
\newcommand{\vV}{\mathcal{V}}

\newcommand{\Supp}{\mathop{\rm Supp}\nolimits}
\newcommand{\Hom}{\mathop{\rm Hom}\nolimits}

\newcommand{\dotimes}{\stackrel{\textbf{L}}{\otimes}}
\newcommand{\dR}{\mathbf{R}}
\newcommand{\dL}{\mathbf{L}}

\newcommand{\Pic}{\mathop{\rm Pic}\nolimits}

\newcommand{\Chow}{\mathop{\rm Chow}\nolimits}

\newcommand{\tX}{\widetilde{X}}

\newcommand{\id}{\textrm{id}}

\newcommand{\ch}{\mathop{\rm ch}\nolimits}
\newcommand{\rk}{\mathop{\rm rk}\nolimits}

\newcommand{\Ext}{\mathop{\rm Ext}\nolimits}
\newcommand{\Spec}{\mathop{\rm Spec}\nolimits}

\newcommand{\Coh}{\mathop{\rm Coh}\nolimits}
\newcommand{\QCoh}{\mathop{\rm QCoh}\nolimits}

\newcommand{\divv}{\mathop{\rm div}\nolimits}

\newcommand{\cneq}{\mathrel{\raise.095ex\hbox{:}\mkern-4.2mu=}}
\newcommand{\eqcn}{\mathrel{=\mkern-4.5mu\raise.095ex\hbox{:}}}
\newcommand{\dit}{\colon}
\newcommand{\Proj}{\mathop{\rm Proj}\nolimits}
\newcommand{\ccdot}{\mathop{\cdot}\nolimits}

\newcommand{\ggcd}{\mathop{\textrm{g.c.d}}\nolimits}

\newcommand{\Bs}{\mathop{\rm Bs}\nolimits}
\newcommand{\Ass}{\mathop{\rm Ass}\nolimits}

\newtheorem{thm}{Theorem}[section]
\newtheorem{prop}[thm]{Proposition}
\newtheorem{lem}[thm]{Lemma}
\newtheorem{defi}[thm]{Definition}
\newtheorem{rmk}[thm]{Remark}
\newtheorem{cor}[thm]{Corollary}
\newtheorem{step}{Step}

\newtheorem{sssstep}{Step}

\newtheorem{ccase}{Case}

\newtheorem{prop-defi}[thm]{Proposition-Definition}

\begin{document}

\maketitle

\begin{abstract}
Let $X$ be a smooth projective variety.
We study a relationship between the derived category of $X$ 
and that of a canonical divisor. 
As an application, we will study Fourier-Mukai transforms
when $\kappa (X)=\dim X -1$. 

\end{abstract}

\section{Introduction}

$\quad$ Let $X$ be a smooth projective variety and $D(X)$ the bounded derived category of 
coherent sheaves on $X$. Recently $D(X)$ draws much attention from many aspects, especially 
mirror symmetry, moduli spaces of stable sheaves, and birational geometry. 
Kontsevich~\cite{Kon} 
conjectures the 
existence of equivalence 
between derived category of $X$ and derived Fukaya category of its mirror. 
In the physical viewpoint, we can not distinguish the mirror pair by observations or experiments, so this gives a motivation for the new concept 
of ``spaces". In this respect, 
the properties which are invariant under Fourier-Mukai transform (i.e. 
categorical invariant) can 
be considered as the essential properties of the ``spaces". 
For example, the Serre functor $S_X =\otimes \omega _X [\dim X]$ 
is such a categorical invariant.

On the other hand there are many works concerning the derived equivalent 
varieties. 
Let $FM(X)$ be a set of isomorphism class of smooth projective varieties which have equivalent derived categories to $X$.
 In~\cite{Mu1}, Mukai showed that 
if $A$ is an abelian variety and $\hat{A}$ is its dual variety, then $\hat{A}$ belongs to $FM(A)$. 
This fact implies that $D(X)$ does not completely determine $X$. 
But if we assume that $K_X$ or $-K_X$ is ample, Bondal and Orlov~\cite{B-O1} showed that $FM(X)$ consists of 
$X$ itself.   
When $X$ is a minimal surface, Bridgeland-Maciocia~\cite{B-M1} described
$FM(X)$, and non-minimal case was treated by Kawamata~\cite{Ka1}. 
In these cases, we can see the following common phenomenon:

``If there exist more information of $K_X$, for example if $\kappa (X, \pm K_X)$ are greater, 
then $FM(X)$ is smaller. "

The main purpose of this paper is to explain why this phenomenon occurs. 
The idea is to extract information concerning Serre functors.  
Here we state the main theorem.  Let $Y\in FM(X)$ and $\Phi \colon D(X)\to D(Y)$ an equivalence of triangulated categories. Let $\pP \in D(X\times Y)$ be
a kernel of $\Phi$. Here the definition of kernel will be given in 
Definition~\ref{ker}. Let $\Psi \colon D(Y)\to D(X)$ be 
a quasi-inverse of $\Phi$, and $\eE \in D(X\times Y)$ be a kernel of $\Psi$.
Then we prove that 
\begin{itemize} 
\item $\Phi$ induces an isomorphism of vector spaces, $H^0 (X,mK_X)\to H^0 (Y,mK_Y)$ for $m\in \mathbb{Z}$. This is also proved in ~\cite{Cal}.
Let $E\in |mK_X|$ corresponds to 
$E^{\dag}\in |mK_Y|$. 
\item $\Phi$ induces a bijection between $\pi _0 (\cap _{i=1} ^n E_i)$ and
$\pi _0 (\cap _{i=1} ^n E_i ^{\dag})$. Here $E_i \in |m_i K_X|$ for $i=1,\cdots ,n$ and $m_i \in \mathbb{Z}$. $n$ and $m_i$ are arbitrary, and $\pi _0$ 
means connected component. Let $C\in 
 \pi _0 (\cap _{i=1} ^n E_i)$ corresponds to $C^{\dag}\in \pi _0 (\cap _{i=1} ^n E_i ^{\dag})$. 
 \end{itemize}
Then the main theorem is the following:
\begin{thm}\label{mt}
Assume that $C$ and $C^{\dag}$ satisfy the following conditions:
\begin{itemize}
\item $C$ and $C^{\dag}$ are complete intersections. 
\item $\pP \dotimes \oO _{C\times Y}$, $\eE \dotimes \oO _{C\times Y}$ 
are sheaves, up to shift. 
\end{itemize}
Then there exists an equivalence of 
triangulated categories
$\Phi _C \colon D(C)\to D(C^{\dag})$ such that the following diagram 
is 2-commutative: 
$$
\begin{CD}
D(X) @>\dL i_C ^{\ast}>>  D(C)  @> i_{C\ast}>>  D(X) \\
@V\Phi VV @V\Phi _C VV @V\Phi VV   \\
D(Y) @>\dL i_{C^{\dag}} ^{\ast}>>  D(C^{\dag})  @> i_{C^{\dag}\ast}>>  D(Y).
\end{CD}
$$
\end{thm}
The assumptions are satisfied if $|m_i K_X|$ are free, $E_i \in |m_i K_X|$
are generic members, and $\pP$ is a sheaf, up to shift.  
The above theorem says that 
``If there are many members in $|mK_X|$, then we can reduce the problem of 
describing $FM(X)$ to lower dimensional case."
As an application, 
we will study Fourier-Mukai transforms when $\kappa (X)=\dim X -1$.
Using this method, we will 
give a generalization of the theorem of Bondal and Orlov~\cite{B-O1}, 
and determine $FM(X)$ when $\dim X=3$ and $\kappa (X)=2$.

In the view point of birational geometry, there are some works concerning 
derived categories and birational geometry. 
For example Bridgeland~\cite{Br1} constructed smooth 3-dimensional flops as a 
moduli space of perverse point sheaves, which are objects in derived
 category. Surprisingly his method gives an equivalence of 
 derived categories under flops simultaneously. 
This result was generalized 
by Chen~\cite{Ch} and Kawamata~\cite{Ka1}. The existence of flops and flips is a very difficult problem in birational geometry,
and Bridgeland's result gives a possibility of treating the problem by moduli theoretic method.

\section{Derived categories and Serre functors}
\subsection*{Notations and conventions}
\begin{itemize}
\item
Throughout this  
paper, we assume all the varieties are defined over $\mathbb{C}$.
\item For smooth projective variety $X$, let $D(X)\cneq D^b (\Coh (X))$, i.e. bounded derived category of coherent sheaves on $X$. The translation functor is written $[1]$, and
the symbol $E[m]$ means the object $E$ shifted to the left by $m$ places.
\item $\omega _X$ means canonical bundle, and $K_X$ means canonical divisor. 
For a Cartier divisors $D$, we write the global section of $\oO _X(D)$
as $H^0 (X,D)$, $|D|$ means linear system, and $\Bs |D|$ is a base locus
 as usual. 
\item For the derived functors, we omit $\dR$ or $\dL$ if the functors we 
want to derive are exact. 
\item For another variety $Y$, we denote 
by $p_i$, projections $p_1 \dit X\times Y \to X$, $p_2 \dit X\times Y \to Y$.
\item For a closed point $x\in X$, $\oO_x$ means a skyscraper 
sheaf supported at 
$x$. 
\end{itemize}

In this section we recall some definitions and properties concerning 
derived categories. 

\begin{defi}\label{ker} For an object
$\mathcal{P}\in D(X\times Y)$, we define a functor $\Phi_{X\to Y} ^{\mathcal{P}}\colon
D(X) \to D(Y)$ by 
$$\Phi_{X\to Y} ^{\mathcal{P}}(E)\cneq \mathbf{R}p_{2\ast}(p_1 ^{\ast}E\stackrel{\mathbf{L}}{\otimes}\mathcal{P}).$$
The object $\mathcal{P}$ is called the kernel of $\Phi_{X\to Y}^{\mathcal{P}}$.
For a morphism $\mu \colon \pP _1 \to \pP _2$ in $D(X\times Y)$, we also denote by 
$\Phi _{X\to Y}^{\mu}$ the natural transform: 
$$\Phi _{X\to Y}^{\mu}\colon \Phi _{X\to Y}^{\pP _1}\longrightarrow \Phi _{X\to Y}^{\pP _2}, $$
induced by $\mu$.
\end{defi}
\label{f}

The functor of the form $\Phi _{X\to Y}^{\pP}$ is called an integral functor. 
If an integral functor gives an equivalence of categories, then it is called a Fourier-Mukai 
transform. The following theorem is fundamental in this paper:

\begin{thm}[Orlov $~\cite{Or1}$]
Let $\Phi \dit D(X) \to D(Y)$ gives an equivalence of $\mathbb{C}$-linear 
triangulated categories. 
Then there exists an object $\mathcal{P} \in D(X\times Y)$ such that $\Phi$ is isomorphic to the functor 
$\Phi_{X\to Y} ^{\mathcal{P}}$. Moreover $\pP$ is uniquely determined up to isomorphism. 
\label{O1}
\end{thm}

Next we introduce the notion of Fourier-Mukai partners.
\begin{defi}
We define $FM(X)$ as the set of isomorphism classes of smooth projective varieties $Y$, which has an 
equivalence of $\mathbb{C}$-linear triangulated categories, $\Phi \dit D(X) \to D(Y)$. 
 If $Y \in FM(X)$, $Y$ is called a Fourier-Mukai partner of $X$.
\end{defi}
 
By Theorem~\ref{O1}, if $Y\in FM(X)$, then
$D(Y)$ is related to $D(X)$ by a Fourier-Mukai 
transform.  
To study the relation between derived categories and canonical divisors, 
the following Serre functor plays an important role. 

\begin{defi}
Let $\mathcal{T}$ be a $\mathbb{C}$-linear triangulated category of 
finite type. An exact equivalence $S\dit \mathcal{T} \to \mathcal{T}$
is called a Serre functor if there exists a bifunctorial isomorphism

$$ \Hom (E,F) \to \Hom (F,S(E))^{\ast}$$
for $E,F \in \mathcal{T}$.
\end{defi}

As in~\cite[Proposition 1.5]{B-O1}, if a Serre functor exists, then it is unique up to canonical isomorphism. If 
$X$ is a smooth projective variety and $\mathcal{T} =D(X)$, then Serre duality implies Serre functor
$S_X$ is given by  $S_X (E) = E \otimes \omega _X [\dim X]$.

\begin{prop-defi} \label{comp}
Let $X$, $Y$, $Z$ be varieties, 
and $p_{ij}$ be projections from 
$X\times Y\times Z$ onto corresponding factors.
Let us take $\fF \in D(X\times Y)$, $\gG\in D(Y\times Z)$. 
We define $\gG \circ \fF \in D(X\times Z)$ as 
$$\gG \circ \fF \cneq \dR p_{13\ast}(p_{12}^{\ast}\fF \dotimes p_{23}^{\ast}\gG ).$$
Then we 
have the isomorphism of functors: $\Phi _{Y\to Z}^{\gG}\circ \Phi _{X\to Y}^{\fF} \cong \Phi _{X\to Z}^{\gG \circ \fF}$, and 
for a morphism
 $\mu \colon 
\fF _1 \to \fF _2$ in $D(X\times Y)$,
the isomorphism of 
natural transforms:
$$\Phi _{Y\to Z}^{\gG}\circ \Phi _{X\to Y}^{\mu} \cong \Phi _{X\to Z}^{\gG \circ \mu}:\Phi _{X\to Z}^{\gG \circ \fF _1} \longrightarrow \Phi _{X\to Z}^{\gG \circ \fF _1}.$$
Moreover the operation $\circ$ is associative, i.e. 
$(\hH \circ \gG)\circ \fF \cong \hH \circ (\gG \circ \fF)$. 
\end{prop-defi}
\begin{proof} The proof of $\Phi _{Y\to Z}^{\gG}\circ \Phi _{X\to Y}^{\fF} \cong \Phi _{X\to Z}^{\gG \circ \fF}$ is seen in several references. 
 For example, see ~\cite[Proposition 2.3]{Ch}. 
 The same proof works for natural transforms, formally replacing 
 $\fF$ by $\mu$. 
 We can check the 
operation $\circ$ is associative by the same method, but we would like to 
give the proof for the lack of references.  
Let $X$, $Y$, $Z$, $W$ be varieties, and take $\fF \in D(X\times Y)$, 
$\gG \in D(Y\times Z)$, and $\hH \in D(Z\times W)$. 
We change the index $p_{ij}$ to $p_{XY}$ etc. 
Let $p_{\ast \ast}$, $q_{\ast \ast}$, $r_{\ast \ast}$, and $s_{\ast \ast}$ 
be projections, given as in the following diagram:\\
\xymatrix{
& {X\times Y\times Z}\ar[dl]_{p_{XY}} \ar[d]^{p_{XZ}} \ar[rd]^{p_{YZ}} & \\
X\times Y & X\times Z & Y\times Z, \\
} \qquad
\xymatrix{
& {X\times Z\times W}\ar[dl]_{q_{XZ}} \ar[d]^{q_{XW}} \ar[rd]^{q_{ZW}} & \\
X\times Z & X\times W & Z\times W, \\
} \\  
\xymatrix{
& {X\times Y\times W}\ar[dl]_{r_{XY}} \ar[d]^{r_{XW}} \ar[rd]^{r_{YW}} & \\
X\times Y & X\times W & Y\times W, \\
} \qquad 
\xymatrix{
& {Y\times Z\times W}\ar[dl]_{s_{YZ}} \ar[d]^{s_{YW}} \ar[rd]^{s_{ZW}} & \\
Y\times Z & Y\times W & Z\times W. \\
}

Let $\pi_{\ast \ast}$ or $\pi _{\ast \ast \ast}$ be projections from 
$X\times Y \times Z \times W$ onto corresponding factors, for example as 
in the following diagram:

\xymatrix{
& & X\times Y \times Z \times W \ar[dll]_{\pi _{XYZ}}
\ar[dl]^{\pi_{XYW}}
\ar[d]^{\pi _{XZW}}
\ar[dr]^{\pi _{YZW}}
& & \\
X\times Y\times Z & X\times Y\times W &  X\times Z\times W
  & Y\times Z\times W \\
 }

\vspace{2mm}

Then $\hH \circ (\gG \circ \fF)$ is calculated as 
\begin{align*}
\hH \circ (\gG \circ \fF) & \cong \dR q_{XW\ast}(q_{XZ}^{\ast}(\gG \circ \fF) 
\dotimes q_{ZW}^{\ast}\hH) \\
& \cong \dR q_{XW\ast}(q_{XZ}^{\ast}\dR p_{XZ\ast}(p_{XY}^{\ast}\fF \dotimes 
p_{YZ}^{\ast}\gG )\dotimes q_{ZW}^{\ast}\hH) \\
& \cong \dR q_{XW\ast}(\dR \pi _{XZW\ast}\pi _{XYZ}^{\ast}
(p_{XY}^{\ast}\fF \dotimes 
p_{YZ}^{\ast}\gG )\dotimes q_{ZW}^{\ast}\hH) \\
& \cong \dR q_{XW\ast}\dR\pi _{XZW\ast}(\pi _{XYZ}^{\ast}
(p_{XY}^{\ast}\fF \dotimes 
p_{YZ}^{\ast}\gG )\dotimes \pi _{XZW}^{\ast}q_{ZW}^{\ast}\hH) \\
& \cong \dR \pi _{XW\ast}(\pi _{XY}^{\ast}\fF \dotimes \pi _{YZ}^{\ast}\gG 
\dotimes \pi _{ZW}^{\ast}\hH). 
\end{align*}

Here the third isomorphism follows from flat base change, and fourth isomorphism from projection formula. Similarly, 
$(\hH \circ \gG)\circ \fF$ is calculated as 
\begin{align*}
(\hH \circ \gG) \circ \fF & \cong 
\dR r_{XW\ast}(r_{YW}^{\ast}(\hH \circ \gG)\dotimes r_{XY}^{\ast}\fF) \\
&\cong \dR r_{XW\ast}(r_{YW}^{\ast}\dR s_{YW\ast}(s_{YZ}^{\ast}\gG \dotimes s_{ZW}^{\ast}\hH )
\dotimes r_{XY}^{\ast}\fF) \\
&\cong \dR r_{XW\ast}(\dR \pi _{XYW\ast}\pi _{YZW}^{\ast}(s_{YZ}^{\ast}\gG \dotimes s_{ZW}^{\ast}\hH )
\dotimes r_{XY}^{\ast}\fF) \\
& \cong \dR \pi _{XW\ast}(\pi _{YZ}^{\ast}\gG 
\dotimes \pi _{ZW}^{\ast}\hH \dotimes \pi _{XY}^{\ast}\fF ).
\end{align*}

Therefore we obtain the isomorphism, $\hH \circ (\gG \circ \fF)\cong
(\hH \circ \gG)\circ \fF$. \end{proof}

Here we give one remark. The category $D(X\times Y)$ is like a category
of functors from $D(X)$ to $D(Y)$. In fact an object $\fF \in D(X\times Y)$ 
corresponds to a functor $\Phi _{X\to Y}^{\fF}$, and a morphism
$\fF \to \gG$ gives a natural transform $\Phi _{X\to Y}^{\fF}\to \Phi _{X\to Y}
^{\gG}$. But as remarked in ~\cite{Cal}, this correspondence is not faithful, i.e. non-trivial morphism $\fF \to \gG$ may induce a trivial natural transform.
Though natural transform is a categorical concept, it is not useful for our 
purpose. So sometimes we use the objects of $D(X\times Y)$ instead of functors, and treat their morphisms as if they are natural transforms.

\section{Moduli spaces of stable sheaves}
In this section, we introduce the notations of the moduli spaces of 
stable sheaves, and recall some properties. These are used for the 
applications of
Theorem~\ref{mt}.  
The details are written in the book~\cite{Hu}.
Let $X$ be a projective scheme and $H$ be a polarization. 
For a non-zero object $E\in \Coh (X)$, its Hilbert polynomial has the following form:
$$\chi (E\otimes H^{\otimes m})=\sum_{i=0}^{d} \frac{\alpha _{i}(E)}{i!}m^i \qquad (\alpha_i (E) \in \mathbb{Z}, d=\dim (\Supp E)).$$ 
We define a rank of $E$ and its reduced Hilbert polynomial by 
$$\rk (E)\cneq \alpha_{d}(E)/\alpha_{d}(\mathcal{O}_X), \qquad 
p(E,H)\cneq \chi (E\otimes H^{\otimes m})/\alpha _{d}(E). $$
Now let us introduce the order on $\mathbb{Q}[m]$ as follows: if $p, p'\in \mathbb{Q}[m]$, 
then $p\le p'$ if and only if $p(m) \le p' (m)$ for sufficiently large $m$. We denote $p<p'$ if $p(m)<p'(m)$ for sufficiently large $m$. 

\begin{defi}
A non-zero object $E\in \Coh (X)$ is said to be $H$-semistable if $E$ is pure, i.e. there exists no subsheaf of dimension lower than $d$, and 
for all subsheaves $F \subsetneq E$, we have $p(F,H)\le p(E,H)$. $E$ is said to be $H$-stable if $E$ is $H$-semistable and 
for all subsheaves $F \subsetneq E$, we have $p(F,H)<p(E,H)$.
\end{defi}

Using the above stability, we can consider the moduli spaces of stable
(semistable) sheaves. Also we can consider the 
relative version of the moduli spaces of such
sheaves, under projective morphism $f\colon X\to S$ and 
$f$-ample divisor $H$.
Let $T$ be a $S$-scheme, and  $p_X \colon X\times _S T \to X $ and $p_T \colon X\times _S T \to T$ be projections. We define a contravariant functor
$\overline{\mathcal{M}}^{H}(X/S)\dit (Sch/S)^{\circ} \to (Sets)$  as follows: 
$$\overline{\mathcal{M}}^{H}(X/S)(T) \cneq \left\{ \begin{array}{l}
 \fF \in \Coh (X\times _S T), \textrm{ which are flat over }T, 
 \textrm{and for }\\
 \textrm{all geometric points }\Spec k(t) \to T,  \mathcal{F}|_{X\times\Spec k(t)} \\
 \textrm{is } p_{X}^{\ast}H|_{X\times\Spec k(t)}\textrm{-semistable.} 
 \end{array}
   \right\}/\sim. $$
Here for $E,E' \in \Coh (X\times _S T)$, the equivalence 
relation $\sim$ is the following:
$$E \sim E' \quad \stackrel{\textrm{def}}{\Leftrightarrow} \quad E \cong E' \otimes p_T ^{\ast}\mathcal{L} \quad \textrm{for some }\mathcal{L} \in \Pic (T).$$ 
Then there exists a projective scheme
$$\overline{M}^{H}(X/S)\to S,$$  which corepresents
$\overline{\mathcal{M}}^{H}(X/S)$. Let $M^{H}(X/S) \subset \overline{M}^{H}(X/S)$ be a subset which 
 corresponds to stable sheaves. It is known that  
$M^{H}(X/S)$ is an open subscheme of $\overline{M}^{H}(X/S)$, for 
example see~\cite{Hu}. 

\begin{defi}
Let $M \subset M^{H}(X/S)$
be an irreducible component. $M$ is called fine if it is projective over $S$ and there exists a universal sheaf on $X\times _S M$. \end{defi} 
The following theorem is due to Mukai~\cite{Mu2}  

\begin{thm}[Mukai $~\cite{Mu2}$]
For $x\in M$, we denote by $E_x$ the corresponding stable sheaf. Then there exists a universal family on $X\times _S M$ if 
$$\emph{g.c.d} \{\chi(E_x \otimes \mathcal{N}) \mid \mathcal{N} \emph{ is a vector bundle on } X \}=1$$ \label{univ}
holds.
\end{thm}

We have the following criteria to find the fine moduli scheme: 

\begin{lem}\label{fine}
If $\emph{g.c.d} \{\chi(E_x \otimes H^{\otimes n}) \mid n\in \mathbb{Z} \}=1$, 
then $M$ is projective over $S$, i.e. there exists no properly semistable 
boundary. Hence $M$ is fine by Theorem~\ref{univ}. \end{lem}
\begin{proof}
Indeed if there exists some  $x\in \overline{M}\setminus M$, { then there exists a subsheaf }$F \subsetneq E_x$ { such that }
$p(F,H)=p(E_x, H)$. { If we take }$n_i, \omega _i \in \mathbb{Z}$ { such that }$\sum \omega _i \cdot \chi (E_x \otimes H^{\otimes n_i})=1$, 
{ then }
$$\sum \omega _i \cdot \chi (F \otimes H^{\otimes n_i})={\alpha _d (F)}/{\alpha _d (E_x) }.$$ { Since the left hand side is an integer and }
$0<{\alpha _d (F)}/{\alpha _d (E_x) }<1$, {we have a contradiction. So by the above theorem }$M$ {is fine.} \end{proof}

Finally we recall the significant result on the moduli spaces of 
stable sheaves and derived categories, established by Bridgeland and
 Maciocia~\cite{B-M2}.
We say that a family of sheaves $\{ \uU _p\}_{p\in M}$ on $X$ is complete if 
the Kodaira-Spencer map
$$T_p M \longrightarrow \Ext ^1 _{X}(\uU _p, \uU _p)$$
is bijective. 

\begin{thm}[Bridgeland-Maciocia ~\cite{B-M2}]\label{int}
Let $X$ be a smooth projective variety of dimension $n$ and $\{ \uU _p\}_{p\in M}$ be a complete family of simple sheaves on $X$ 
parameterized by an irreducible projective scheme $M$ of dimension $n$. Suppose that 
$\Hom _X (\uU _{p_1}, \uU _{p_2})=0$ for $p_i \in M$, $p_1 \neq p_2$ and the set $$\Gamma (\uU )\cneq \{ (p_1, p_2)\in M\times M \mid \Ext _X ^i (\uU _{p_1}, \uU _{p_2})\neq 0 \quad \textrm{for some }i\in \mathbb{Z}\}$$
has $\dim \Gamma (\uU)\le n+1$. Suppose also that $\uU _p \otimes \omega _X \cong \uU _p$ for all $p\in M$. Then $M$ is a nonsingular projective 
variety and $\Phi _{M\to X}^{\uU}\colon D(M)\to D(X)$ is an equivalence.
\end{thm}

\section{Correspondences of canonical divisors}
In this section we fix two smooth projective varieties $X$ and $Y$, such that 
$Y\in FM(X)$. The purpose of this section is to establish the relation between 
the canonical divisors of $X$ and $Y$, and state our main theorem. 
We fix the following notation:

\begin{itemize}
\item $\Phi\colon D(X)\to D(Y)$ gives an equivalence, and $\pP \in D(X\times Y)$ is a kernel of $\Phi$. 
\item $\Psi \colon D(Y)\to D(X)$ is a quasi-inverse of $\Phi$, and 
$\eE \in D(X\times Y)$ is a kernel of $\Psi$. 
\item $S_X \cneq \otimes \omega _X [\dim X] \colon D(X)\to D(X)$ is a Serre 
functor of $D(X)$. 
\end{itemize}
Since Serre functor is categorical, we have the isomorphism of 
functors,
$$\tau \colon \Phi \circ S_X  \stackrel{\sim}{\longrightarrow} S_Y  \circ \Phi .$$
Note that the kernel of left hand side is ${\pP \otimes p_1 ^{\ast}\omega _X}[\dim X]$, and right hand side is $\pP \otimes p_2 ^{\ast}\omega _Y [\dim Y]$. So by Theorem~\ref{O1},
 we have the isomorphism, 
$$\rho \colon \pP \otimes p_1 ^{\ast}\omega _X [\dim X]\stackrel{\sim}{\longrightarrow}\pP \otimes p_2 ^{\ast}\omega _Y [\dim Y].$$
Therefore $\dim X =\dim Y$, and there exist isomorphism for all $m\in \mathbb{Z}$, 
$$\rho _m \colon \pP \otimes p_1 ^{\ast}\omega _X ^{\otimes m}
\stackrel{\sim}{\longrightarrow}\pP \otimes p_2 ^{\ast}\omega _Y ^{\otimes m}.$$Therefore we can see the following proposition:
\begin{prop}\label{md}
$\{ \rho _{m} \} _{m\in \mathbb{Z}}$ induce the isomorphism of graded $\mathbb{C}$-algebras:
$$\{ \rho _{m}' \}\colon \bigoplus _{m\in \mathbb{Z}}\Hom _{X\times Y}(\pP, \pP \otimes p_1 ^{\ast}\omega _X ^{\otimes m}) \stackrel{\sim}{\longrightarrow} 
\bigoplus _{m\in \mathbb{Z}}\Hom _{X\times Y}(\pP, \pP \otimes p_2 ^{\ast}\omega _Y ^{\otimes m}).$$
\end{prop}
\begin{proof} Clear by the above argument. \end{proof}
Next we will compare the vector spaces $H^0 (X,mK_X)$ and $\Hom _{X\times Y}(\pP, \pP \otimes p_1 ^{\ast}\omega _X ^{\otimes m})$.
Since $\Phi$ gives an identification of categories, $\Phi$ must give the bijection between functors $D(X)\to D(X)$ and functors $D(X)\to D(Y)$. 
In this respect, the following lemma is obvious:
\begin{lem}\label{ob}
The following functor,
$$\pP \circ \colon D(X\times X) \ni a \mapsto \pP \circ a \in D(X\times Y)$$ 
gives equivalence. 
\end{lem}
\begin{proof} Let $\Psi$ be a quasi-inverse of $\Phi$, and $\eE \in D(X\times Y)$ be a kernel of $\Psi$. 
Let $\Delta _X \subset X\times X$ and $\Delta _Y \subset Y\times Y$ be diagonals. Note that the operations $\oO _{\Delta _X}\circ$, $\oO _{\Delta _Y}\circ$ induce identities. 
Since $\eE \circ \pP \cong \oO _{\Delta _X}$, 
$\pP \circ \eE \cong \oO _{\Delta _Y}$, the following functor: 
$$\eE \circ \colon D(X\times Y) \ni b \mapsto \eE \circ b \in D(X\times X)$$
gives a quasi-inverse by Proposition-Definition~\ref{comp}. \end{proof}
As in the same way, we have equivalence of categories:
$$\circ \pP \colon D(Y\times Y)\ni a \mapsto a\circ \pP \in D(X\times Y).$$
We have the following lemma:
\begin{lem}\label{co}
The following diagrams are 2-commutative:
$$\begin{CD}
D (X\times X) @>{\pP \circ}>> D (X\times Y) \\
@A{\Delta _{\ast}}AA   @AA{p_1 ^{\ast}(\ast)\dotimes \pP}A \\
D (X) @= D (X),
\end{CD} \qquad \quad
\begin{CD}
D (Y\times Y) @>{\circ\pP}>> D (X\times Y) \\
@A{\Delta _{\ast}}AA   @AA{p_2 ^{\ast}(\ast)\dotimes \pP}A \\
D (Y) @= D (Y).
\end{CD}$$
Here $\Delta$ means diagonal embedding. 
\end{lem}
\begin{proof} 
Let us check the left diagram commutes. 
Let $p_{ij}$ be projections from $X\times X\times Y$ onto 
corresponding factors.
Take $a\in D (X)$. Then 
\begin{align*}
\pP \circ (\Delta _{\ast}a) &\cong \dR p_{13\ast}\left( p_{12}^{\ast}\Delta _{\ast}a \dotimes p_{23}^{\ast}\pP \right) \\
&\cong \dR p_{13\ast}\left( (\Delta \times \id _Y)_{\ast}p_1 ^{\ast}a \dotimes p_{23}^{\ast}\pP \right) \\
&\cong \dR p_{13\ast}(\Delta \times \id _Y)_{\ast}\left( p_1 ^{\ast}a \dotimes (\Delta \times \id _Y)^{\ast}p_{23}^{\ast}\pP \right) \\
&\cong p_1 ^{\ast}a \dotimes \pP.
\end{align*}
The second isomorphism follows from 
flat base change of the diagram below 
$$\begin{CD}
X\times Y @>{\Delta \times \id _Y}>> X\times X\times Y \\
@V{p_1}VV @VV{p_{12}}V \\
X @>{\Delta }>> X\times X, 
\end{CD}$$
and the third isomorphism follows from projection formula. \end{proof} 
As the immediate corollary, we have 
\begin{cor}\label{cano}
$\Phi$ induces the isomorphism of graded $\mathbb{C}$-algebras:
$$\{ \phi _{m} \} _{m\in \mathbb{Z}} \colon 
\bigoplus _{m\in \mathbb{Z}}H^0 (X,mK_X) \stackrel{\sim}{\longrightarrow}
\bigoplus _{m\in \mathbb{Z}}H^0 (Y,mK_Y).$$
\end{cor}
\begin{proof}
By Lemma~\ref{co}, we have the isomorphism of graded $\mathbb{C}$-algebras:
\begin{align*}\bigoplus _{m\in \mathbb{Z}}\Hom _{X\times X}
(\Delta _{\ast}\oO _X, 
\Delta _{\ast}\omega _X ^{\otimes m}) \stackrel{\pP \circ}{\longrightarrow} 
\bigoplus _{m\in \mathbb{Z}}\Hom _{X\times Y}
(\pP, \pP \otimes p_1 ^{\ast}
\omega _X ^{\otimes m}), \\
\bigoplus _{m\in \mathbb{Z}}\Hom _{Y\times Y}
( \Delta _{\ast}\oO _Y ,\Delta _{\ast}\omega _Y ^{\otimes m}) \stackrel{\circ \pP}{\longrightarrow} 
 \bigoplus _{m\in \mathbb{Z}}
\Hom _{X\times Y}(\pP, \pP \otimes p_2 ^{\ast}
\omega _Y ^{\otimes m}).
\end{align*}
Since $H^0 (X,mK_X)=\Hom _{X\times X}(\Delta _{\ast}\oO _X, 
\Delta _{\ast}\omega _X ^{\otimes m})$, combining $\rho ' _{m}$ given in 
Proposition~\ref{md}, we obtain the corollary.
\end{proof}
Now let us interpret the isomorphism $\phi _{m} \colon H^0 (X,mK_X) \to H^0 (Y,mK_Y)$
categorically. 
Take $\sigma \in H^0 (X,mK_X)$ and $\sigma ^{\dag}\cneq \phi _{m} (\sigma)\in H^0 (Y,mK_Y)$. Let $d\cneq \dim X =\dim Y$. 
Then we can think $\sigma$ and $\sigma ^{\dag}$ as natural transforms, 
$$\sigma \colon \id _X \longrightarrow
 S_X ^m [-md], \quad \sigma ^{\dag}\colon 
\id _Y \longrightarrow S_Y ^m [-md]$$
Here $S_X ^m [-md]$ is a $m$-times composition of the shifted Serre functor, 
$S_X [-d]=\otimes \omega _X$. Let 
$$\tau _m \colon \Phi \circ S_X ^m [-md] \stackrel{\sim}{\longrightarrow}
S_Y ^m [-md] \circ \Phi$$
be the isomorphism of functors, induced by $\tau \colon \Phi \circ S_X \stackrel{\sim}{\longrightarrow} S_Y \circ \Phi$ naturally. 

\begin{lem}
$\sigma ^{\dag}$ is equal to the following composition:
$$\emph{id} _Y=\Phi \circ \emph{id} _X \circ \Phi ^{-1}\stackrel{\emph{id} \circ \sigma \circ \emph{id}}{\longrightarrow}\Phi \circ S_X ^m [-dm]\circ \Phi ^{-1} 
\stackrel{\tau _m \circ \emph{id}}{\longrightarrow}S_Y ^m [-dm]\circ \Phi \circ \Phi ^{-1}=S_Y ^m [-dm].$$
\end{lem}
\begin{proof}
This follows from Proposition-Definition~\ref{comp} for natural transforms, 
and the construction of 
$\sigma ^{\dag}$. \end{proof}

Let $E\cneq \divv (\sigma)\in |mK_X|$, $E^{\dag}\cneq \divv (\sigma ^{\dag})
\in |mK_Y|$. 
For the closed subscheme $Z\hookrightarrow X$, we define the full subcategory $D_Z (X)\subset D(X)$ as follows:
$$D_Z (X)\cneq \{ a\in D(X) \mid \Supp a \cneq \cup \Supp H^i(a)\subset Z \}.$$
We can observe the following:
\begin{lem}\label{support}
$\Phi$ takes $D_E (X)$ to $D_{E^{\dag}}(Y)$.\end{lem}
\begin{proof}
Take $a\in \Coh (X) \cap D_E (X)$. Let $\sigma ^l \colon \id _X 
\to S_X ^{lmd}[-lmd]=\otimes \omega _X ^{\otimes lm}$ be 
$l$-times composition of 
$\sigma$. Then 
$$\sigma ^l (a)\colon a \to a \otimes \omega _X ^{\otimes lm}$$
are zero-maps for sufficiently large $l$. Then by the above categorical interpretation of $\sigma ^{\dag}$, we have 
$$(\sigma ^{\dag})^l(\Phi (a))\colon \Phi (a) \to \Phi (a) \otimes \omega _Y ^{\otimes lm}$$
are also zero-maps. Since  $(\sigma ^{\dag})^l$ is a natural transform, locally multiplying the defining equation of $lE^{\dag}$, we have
$\Supp \Phi (a) \subset E^{\dag}$. Since
$D_E (X)$ is generated by $\Coh (X) \cap D_E (X)$, the lemma follows. \end{proof}

For the sake of applications, it is convenient to generalize the above lemma to the intersections of canonical divisors. 
\begin{cor}\label{support2}
Take $E_i \in |m_i K_X|$ and their corresponding divisors $E_i ^{\dag}\in |m_i K_Y|$ for $i=1,2 \cdots ,n$. 
There exists a one-to-one correspondence, 
$$\pi _0 (\cap _{i=1}^n E_i ) \ni C \longmapsto C^{\dag}\in 
\pi _0 (\cap _{i=1}^n E_i ^{\dag} ), $$
such that $\Phi$ takes $D_C (X)$ to $D_{C^{\dag}}(Y)$. Here $\pi _0$ means 
connected component. 
\end{cor}
\begin{proof}
Lemma~\ref{support} shows that $\Phi$ takes $D_{\cap E_i}(X)$ to $D_{\cap E_i ^{\dag}}(Y)$. 
Take a connected component $C\subset \bigcap _{i=1}^n E_i$. Since 
$$\Hom _Y(\Phi(\mathcal{O}_{C_{\textrm{red}}}), \Phi(\mathcal{O}_{C_{\textrm{red}}}))
=\Hom _X(\mathcal{O}_{C_{\textrm{red}}},\mathcal{O}_{C_{\textrm{red}}})=\mathbb{C},$$ 
$\Supp \Phi(\mathcal{O}_{C_{\textrm{red}}})$ is connected. Therefore there exists a unique connected component
$C^{\dag}\subset \bigcap _{i=1}^n E_i ^{\dag}$ such that $\Supp \Phi(\mathcal{O}_{C_{\textrm{red}}})\subset C^{\dag}$. 
We show $\Phi$ takes $D_C (X)$ to 
$D_{C^{\dag}}(Y)$. 
It suffices to show $\Phi$ takes $\Coh (\oO _{C})$ to $D_{C^{\dag}}(Y)$. Take a closed point $x\in C$. Then $\Supp (\Phi (\oO _x))$ is connected by the same reason. Since there exists a non-trivial morphism $\oO _{C_{\textrm{red}}}\to \oO _x$, we have $\Phi (\oO _x)\in D_{C^{\dag}}(Y)$. Let us take a simple $\oO _C$-module $\fF$. Then since $\Supp (\Phi( \fF ))$ is connected and there exists a 
non-trivial morphism $\fF \to \oO _x$ for some closed point $x\in C$, we have 
$\Phi (\fF)\in D_{C^{\dag}}(Y)$. The lemma follows by taking 
Harder-Narasimhan filtrations. 
\end{proof}

\vspace{5mm}

Unfortunately the natural functor $D(C)\to D_C (X)$ does not give an equivalence. (In general, latter has bigger $\Ext$-groups.) 
However the existence of equivalence between $D_C (X)$ and $D_{C^{\dag}}(Y)$ leads us to the speculation that $D(C)$ and $D(C^{\dag})$ may be equivalent.  
If $D(C)$ and $D(C^{\dag})$ are equivalent, then the relation between 
$C$ and $C^{\dag}$ will give us information of the relation between $X$ and 
$Y$. One of the purpose of this paper is to claim this speculation is true, under some technical conditions. 
We assume the following conditions on $C$, $C^{\dag}$, and $\pP, \eE \in D(X\times Y)$. Recall that $\pP, \eE$ are kernels of $\Phi$ and $\Phi ^{-1}$. 
\begin{itemize}
\item $C$ and $C^{\dag}$ are complete intersections. 
\item $\pP \dotimes \oO _{C\times Y}$ and $\eE \dotimes \oO _{C\times Y}$ 
are sheaves, up to shift. 
\end{itemize}
These conditions are satisfied, for example, $|m_i K_X|$ are free and $E_i$ are generic members, and $\pP$ is a sheaf. 
Now we can state our main theorem. 
\begin{thm} \label{MT}
Under the above conditions, there exists equivalence 
$\Phi _C \colon D(C) \to D(C^{\dag})$ such that the following diagram is 2-commutative:
$$
\begin{CD}
D(X) @>\dL i_C ^{\ast}>>  D(C)  @> i_{C\ast}>>  D(X) \\
@V\Phi VV @V\Phi _C VV @V\Phi VV   \\
D(Y) @>\dL i_{C^{\dag}} ^{\ast}>>  D(C^{\dag})  @> i_{C^{\dag}\ast}>>  D(Y).
\end{CD}
$$
Here $i_C$, $i_{C^{\dag}}$ are inclusions of $C$, $C^{\dag}$ into $X$ and $Y$ respectively.
\end{thm}

\section{Proof of Theorem~\ref{MT}}

In this section, we give the proof of Theorem~\ref{MT}. 
We use the notation of the 
previous section. 
At first, we explain the plan of the proof.  
We will divide the proof into 4-Steps. In Step 1 and 2, we will show there exists an isomorphism, $\pP \dotimes \oO _{C\times Y}\cong \pP \dotimes \oO _{X\times C^{\dag}}$. Using this and the assumptions, we will find a $\pP _C \in D(C\times C^{\dag})$, and 
construct a functor $\Phi _C \colon D(C) \to D(C^{\dag})$. In Step 4 and 5, we will show $\Phi _C$ gives the desired equivalence.

\begin{sssstep}\label{1}
There exists an isomorphism $\pP \dotimes \oO _{E_i \times Y}\cong \pP 
\dotimes \oO _{X\times E_i ^{\dag}}$.  
\end{sssstep}
\begin{proof}
We omit the index $i$, and write $E_i$ as $E$, etc. 
We have the following exact sequences:
$$ \begin{array}{ccccccccc} 0 &{\longrightarrow}& \mathcal{O}_X 
&\stackrel{\sigma}{\longrightarrow}& \omega _X ^{\otimes m} &\longrightarrow& \mathcal{O}_{E}\otimes \omega _X ^{\otimes m} &\longrightarrow& 0 \\
0 &{\longrightarrow}& \mathcal{O}_Y 
&\stackrel{\sigma ^{\dag}}{\longrightarrow}& \omega _Y ^{\otimes m} &\longrightarrow& \mathcal{O}_{E^{\dag}}\otimes \omega _Y ^{\otimes m} &\longrightarrow& 0 .\end{array}  $$
Applying $p_1 ^{\ast}(\ast)\dotimes \pP$ and $p_2 ^{\ast}(\ast)\dotimes \pP$ 
respectively, we obtain the distinguished triangles:
$$ \begin{array}{ccccccccc} \pP 
&\stackrel{\id \otimes p_1 ^{\ast} \sigma}{\longrightarrow}& \pP \otimes p_1 ^{\ast}\omega _X ^{\otimes m} &\longrightarrow& \pP \dotimes \mathcal{O}_{E\times Y}\otimes p_1 ^{\ast}\omega _X ^{\otimes m} &\longrightarrow& \pP [1] \\
\pP 
&\stackrel{\id \otimes p_2 ^{\ast} \sigma ^{\dag}}{\longrightarrow}& \pP \otimes p_2 ^{\ast}\omega _Y ^{\otimes m} &\longrightarrow& \pP \dotimes \mathcal{O}_{X\times E^{\dag}}\otimes p_2 ^{\ast}\omega _Y ^{\otimes m} &\longrightarrow& \pP [1]. 
\end{array}  $$
On the other hand, by Lemma~\ref{co} and the definition of $\phi _{m}$
 given in Corollary~\ref{cano}, we obtain the following commutative diagram:
$$\begin{CD}
\pP @>\id \otimes p_1 ^{\ast} \sigma >> \pP \otimes p_1 ^{\ast}\omega _X ^{\otimes m} @>>> \pP \dotimes \mathcal{O}_{E\times Y}\otimes p_1 ^{\ast}\omega _X ^{\otimes m} \\
@| @VV {\rho _m} V @.\\
\pP @>\id \otimes p_2 ^{\ast} \sigma ^{\dag}>> \pP \otimes p_2 ^{\ast}\omega _Y ^{\otimes m} @>>> \pP \dotimes \mathcal{O}_{X\times E^{\dag}}\otimes p_2 ^{\ast}\omega _Y ^{\otimes m}.
\end{CD}$$
Here $\rho _m$ is an isomorphism constructed in the previous section. 
Therefore there exists a (not necessary unique) isomorphism,
$$\mathcal{P}\stackrel{\mathbf{L}}{\otimes}
\mathcal{O}_{E \times Y}\otimes p_1 ^{\ast}\omega _X ^{\otimes m}
\cong\mathcal{P}\stackrel{\mathbf{L}}{\otimes}
\mathcal{O}_{X\times E ^{\dag}}\otimes p_2 ^{\ast}\omega _Y ^{\otimes m}.$$ 
Since $\pP \otimes p_1 ^{\ast}\omega _X ^{\otimes m} \cong \pP \otimes p_2 ^{\ast}\omega _Y ^{\otimes m}$, we have an isomorphism, $\pP \dotimes \oO _{E\times Y}\cong \pP \otimes \oO _{X\times E^{\dag}}$.

\begin{sssstep}\label{2}
There exists an isomorphism, 
$$\pP \dotimes \oO _{C\times Y}\cong \pP \dotimes \oO _{X\times C^{\dag}}.$$
\end{sssstep}
\begin{proof}
By using the isomorphism of Step~\ref{1} 
$n$-times, we can get the isomorphism: 
$$\pP \dotimes \left(\bigotimes _{1 \le i \le n}^{\mathbf{L}}\oO _{E_i \times Y}\right)\cong \pP \dotimes 
\left(\bigotimes _{1 \le i \le n} ^{\mathbf{L}}\mathcal{O}_{X\times E_i ^{\dag}}\right).
$$
On the other hand, we have 
$$\bigotimes _{1 \le i \le n} ^{\mathbf{L}}\mathcal{O}_{E_i \times Y}=\bigoplus _{C\in \pi _0 (\bigcap _{i=1} ^{n} E_i)}p_1 ^{\ast}\mathcal{A}_C , \quad
\bigotimes _{1 \le i \le n} ^{\mathbf{L}}\mathcal{O}_{X\times E_i ^{\dag}}=
\bigoplus _{C' \in \pi _0 (\bigcap _{i=1} ^{n} E_i ^{\dag})}p _2 ^{\ast}\mathcal{B}_{C'} ,$$
for some $\aA _C \in D_C(X)$, $\bB _{C'}\in D_{C'}(Y)$. Therefore we have the following isomorphism:
$$  
\bigoplus _{C\in \pi _0 (\bigcap _{i=1} ^{n} E_i)}\mathcal{P}\stackrel{\mathbf{L}}{\otimes}p_1 ^{\ast}\mathcal{A}_C 
\cong  
\bigoplus _{C' \in \pi _0 (\bigcap _{i=1} ^{n} E_i ^{\dag})}\mathcal{P}\stackrel{\mathbf{L}}{\otimes}p_2 ^{\ast}\mathcal{B}_{C'} .
$$
Now we have the following lemma:
\begin{lem}
$\pP \dotimes p _1 ^{\ast}\aA _C$, 
$\pP \dotimes p_2 ^{\ast}\bB _{C^{\dag}}$
are supported on $C\times C^{\dag}$. 
\end{lem}
\begin{proof}
We show $\pP \dotimes p_1 ^{\ast}\aA _C$ is supported on $C\times C^{\dag}$. 
The rest case follows similarly. We can write, 
$$\pP \dotimes p_1 ^{\ast}\aA _C \cong \bigoplus _{C' \in \pi _0 (\bigcap _{i=1} ^{n} E_i ^{\dag})}\rR _{C'}, $$
where $\rR _{C'}$ is supported on $C\times C'$. Take $C' \neq C \in \pi _0 (\cap E_i)$ and assume $\rR _{C'}$ is not zero. 
Let us take a sufficiently ample 
line bundle $\lL$ on $X$. Since $\Phi (\aA _C \otimes \lL) \in D_{C^{\dag}}(Y)$, we have $\dR p_{2\ast}(\rR _{C'}\otimes p_1 ^{\ast}\lL)=0$. On the 
other hand, if $\lL$ is sufficiently ample and $H^q (\rR _{C'})\neq 0$, then
$p_{2\ast}(H^{q}(\rR _{C'})\otimes p_1 ^{\ast}\lL)\neq 0$ and 
$R ^p p_{2\ast}(H^{q}(\rR _{C'})\otimes p_1 ^{\ast}\lL)=0$ for $p>0$. 
Since there exists a following spectral sequence:
$$E_2 ^{p,q}=R ^p p_{2\ast}(H^{q}(\rR _{C'})\otimes p_1 ^{\ast}\lL)
\Rightarrow \dR ^{p+q}p_{2\ast}(\rR _{C'}\otimes p_1 ^{\ast}\lL),$$ we have
$\dR p_{2\ast}(\rR _{C'}\otimes p_1 ^{\ast}\lL)\neq 0$. But this is a 
contradiction. 
 \end{proof}
By the lemma above we have $\pP \dotimes p_1 ^{\ast}\aA _C \cong \pP \dotimes p_2 ^{\ast}\bB _{C^{\dag}}$. 
Since we have assumed $C$ and $C^{\dag}$ are complete intersections, we have $\aA _C =\oO _C$, $\bB _{C^{\dag}}=\oO _{C^{\dag}}$ in our case. Combining these, we have the desired isomorphism:
$$\mathcal{P}\stackrel{\mathbf{L}}{\otimes}\mathcal{O}_{C\times Y} \cong \mathcal{P}\stackrel{\mathbf{L}}{\otimes}\mathcal{O}_{X\times C^{\dag}}.$$
\end{proof}

By the assumptions, the object
$\pP \dotimes \oO _{C\times Y}\cong \pP \dotimes \oO _{X\times C^{\dag}}$
is a sheaf, up to shift. 
This sheaf is $\oO _{C\times Y}$-module and also 
$\oO _{X\times C^{\dag}}$-module. 
Hence this object is a $\oO _{C\times C^{\dag}}$-module, 
so there exists an object $\pP _{C}\in D(C\times C^{\dag})$, such that
$$\pP \dotimes \oO _{C\times Y}\cong \pP \dotimes \oO _{X\times C^{\dag}}
\cong i_{C\times C^{\dag}\ast}\pP _C.$$
Let
$\Phi _C \cneq \Phi_{C\to C^{\dag}}^{\mathcal{P}_C}
\colon D(C)\to D(C^{\dag})$. 
In what follows, we don't use the fact these are sheaves up to shift, and 
show that $\Phi _C$ gives a desired equivalence.

\begin{sssstep}\label{4}
In the diagram of Theorem~\ref{MT} we have the following isomorphisms of functors:
\begin{align*} \Phi _C \circ \mathbf{L}i_C ^{\ast} \cong \Phi_{X\to C^{\dag}}^{(i_C \times \emph{id}_{C^{\dag}})_{\ast}\mathcal{P}_C}, &\qquad& 
\mathbf{L}i_{C^{\dag}}^{\ast} \circ \Phi \cong \Phi_{X\to C^{\dag}}^{\mathbf{L}(\emph{id}_X \times i_{C^{\dag}})^{\ast} \mathcal{P}}, \\
i_{C^{\dag}\ast}\circ \Phi _C \cong \Phi_{C\to Y}^{(\emph{id}_C \times i_{C^{\dag}})_{\ast}\mathcal{P}_C}, &\qquad&
\Phi \circ i_{C\ast} \cong \Phi_{C\to Y}^{\mathbf{L}(i_C \times  \emph{id}_Y)^{\ast}\mathcal{P}}.\end{align*}
$($ See the following diagram. $)$

$$\xymatrix{
& X\times C^{\dag}\ar[rd]^{\emph{id}_X \times i_{C^{\dag}}} & \\
C\times C^{\dag}\ar[ru]^{i_C \times \emph{id}_{C^{\dag}}} \ar[rr]^{i_{C\times C^{\dag}}}  \ar[rd]_{\emph{id}_C \times i_{C^{\dag}}}& & X\times Y \\
& C\times Y \ar[ru]_{i_C \times \emph{id}_Y} & \\
}$$

\end{sssstep}
\begin{proof} 
 Let us calculate $\Phi _C \circ \mathbf{L}i_C ^{\ast}$ by using Proposition-Definition~\ref{comp}. The rest formulas follow similarly. 
 Let $q_{12}\dit X\times C\times C^{\dag} \to X\times C$,   
$q_{23}\dit X\times C\times C^{\dag} \to C\times C^{\dag}$,
$q_{13}\dit X\times C\times C^{\dag} \to X\times C^{\dag}$ be projections. 
Let $\Gamma _C \subset X\times C$ be the graph of the inclusion $i_C$. 
Let $j$ be the inclusion of 
 $\Gamma _C \times C^{\dag}$ into $X\times C\times C^{\dag}$.
Since $\mathbf{L}i_C ^{\ast}=\Phi_{X \to C} ^{\mathcal{O}_{\Gamma _C}}$, we can compute the kernel of 
 $\Phi _C \circ \mathbf{L}i_C ^{\ast}$ as follows: 
 \begin{align*}
 \mathbf{R}q_{13 \ast}(q_{12}^{\ast}\mathcal{O}_{\Gamma _C}
 \stackrel{\mathbf{L}}
 {\otimes}q_{23}^{\ast}\mathcal{P}_C) 
 &\cong \mathbf{R}q_{13 \ast}(\mathcal{O}_{\Gamma _C \times C^{\dag}}\stackrel{\mathbf{L}}
 {\otimes}q_{23}^{\ast}\mathcal{P}_C) \\
 &\cong \mathbf{R}q_{13 \ast}j_{\ast}\mathbf{L}j^{\ast}\mathbf{L}q_{23}^{\ast}\mathcal{P}_C \\
 &\cong (i_C \times \textrm{id}_{C^{\dag}})_{\ast} 
 \mathbf{R}q_{23 \ast}j_{\ast}\mathbf{L}j^{\ast}\mathbf{L}q_{23}^{\ast}
 \mathcal{P}_C \\
 &\cong (i_C \times \textrm{id}_{C^{\dag}})_{\ast}\dR (q_{23}\circ j)_{\ast}\dL (q_{23}\circ j)^{\ast}\pP _C \\
 &\cong (i_C \times \textrm{id}_{C^{\dag}})_{\ast}\mathcal{P}_C .
 \end{align*} 
Here the third isomorphism follows from $q_{13}\circ j =(i_C \times \textrm{id}_{C^{\dag}})\circ q_{23}\circ j$ and the
last isomorphism follows since $q_{23}\circ j$ is identity.   \end{proof}

\begin{sssstep} $\Phi _C$ gives a desired equivalence. \label{5}
\end{sssstep}
By Step~\ref{4}, to prove the diagram of Theorem~\ref{MT} commutes, we only have to check the followings hold:
$$(i_C \times \textrm{id}_{C^{\dag}})_{\ast}\mathcal{P}_C \cong \mathbf{L}(\textrm{id}_X \times i_{C^{\dag}})^{\ast} \mathcal{P},
\qquad 
(\textrm{id}_C \times i_{C^{\dag}})_{\ast}\mathcal{P}_C \cong 
\mathbf{L}(i_C \times  \textrm{id}_Y)^{\ast}\mathcal{P}.$$
There exists a following morphism:
$$\mathcal{P} \to \mathcal{P}\stackrel{\mathbf{L}}{\otimes}\mathcal{O}_{X\times C^{\dag}}\cong i_{C\times C^{\dag} \ast}\mathcal{P}_C
=(\textrm{id}_X \times i_{C^{\dag}})_{\ast}(i_C \times \textrm{id}_{C^{\dag}})_{\ast}\mathcal{P}_C .$$
Taking its adjoint, we have a morphism $\mathbf{L}(\textrm{id}_X \times i_{C^{\dag}})^{\ast} \mathcal{P} \to 
(i_C \times \textrm{id}_{C^{\dag}})_{\ast}\mathcal{P}_C $. Let us take its distinguished triangle: 
$$\qQ\to \mathbf{L}(\textrm{id}_X \times i_{C^{\dag}})^{\ast} \mathcal{P} \to 
(i_C \times \textrm{id}_{C^{\dag}})_{\ast}\mathcal{P}_C \to \qQ[1].$$
By applying $(\textrm{id}_X \times i_{C^{\dag}})_{\ast}$, we get the distinguished triangle,
$$(\textrm{id}_X \times i_{C^{\dag}})_{\ast}\qQ \to 
\mathcal{P}\stackrel{\mathbf{L}}{\otimes}\mathcal{O}_{X\times C^{\dag}}\stackrel{\cong}{\to} i_{C\times C^{\dag} \ast}\mathcal{P}_C
\to (\textrm{id}_X \times i_{C^{\dag}})_{\ast}\qQ[1].$$
So, we have $(\textrm{id}_X \times i_{C^{\dag}})_{\ast}\qQ=0$. Therefore $\qQ=0$ and the morphism 
$\mathbf{L}(\textrm{id}_X \times i_{C^{\dag}})^{\ast} \mathcal{P} {\to} 
(i_C \times \textrm{id}_{C^{\dag}})_{\ast}\mathcal{P}_C $ is an isomorphism.  
We can prove the isomorphism, $(\textrm{id}_C \times i_{C^{\dag}})_{\ast}\mathcal{P}_C \cong 
\mathbf{L}(i_C \times  \textrm{id}_Y)^{\ast}\mathcal{P}$ similarly. 

Finally, we prove $\Phi _C$ gives an equivalence. 
Let us define $\Psi _C \colon D(C^{\dag})\to D(C)$ as
in the same way of $\Phi _C$, from $\Psi = \Phi ^{-1}$. Then 
the following diagram commutes: 
$$
\begin{CD}
D(X) @>\dL i_C ^{\ast}>>  D(C)  @>i_{C\ast}>>  D(X) \\
@V\Phi VV @V\Phi _C VV @V\Phi VV   \\
D(Y) @>\dL i_{C^{\dag}} ^{\ast}>>  D(C^{\dag})  @>i_{C^{\dag}\ast}>>  D(Y) \\
@V\Psi VV @V\Psi _{C} VV @V\Psi VV   \\
D(X) @>\dL i_C ^{\ast}>>  D(C)  @>i_{C\ast}>>  D(X).
\end{CD}
$$
Take a closed point $x \in C$. 
Then by the diagram above, 
 $i_{C\ast}\circ \Psi _C\circ \Phi _C (\mathcal{O}_x)\cong i_{C\ast}(\mathcal{O}_x)$, so $
\Psi _C\circ \Phi _C (\mathcal{O}_x) \cong \mathcal{O}_x$. Then, by ~\cite[Lemma 4.3]{Br2}, 
kernel of $\Psi _C\circ \Phi _C$ is a sheaf on $C\times C$, therefore it must be a line bundle on its 
diagonal. Hence $\Psi _C\circ \Phi _C \cong \otimes \mathcal{L}_C $ for some line bundle $\mathcal{L}_C$ on $C$. 
But, again by the diagram above, we have $\Psi _C\circ \Phi _C (\mathcal{O}_C)\cong \mathcal{O}_C$. This implies 
$\mathcal{L}_C \cong \mathcal {O}_C$ and $\Psi _C \circ \Phi _C \cong \textrm{id}$. 
Similarly, $\Phi _C\circ \Psi _C\cong \textrm{id}$. Therefore $\Phi _C$ is an equivalence and the proof of 
Theorem~\ref{MT} is completed.  \end{proof}

\begin{rmk}
The conditions of kernels are required to find the 
object $\pP _C$ which satisfies 
$$\pP \dotimes \oO _{C\times Y}\cong 
\pP \dotimes \oO _{X\times C^{\dag}}\cong i_{C\times C^{\dag}\ast}\pP _{C}.$$
In fact, if we can find such a $\pP _C$, then our theorem holds by using 
$\pP _{C}$. In Step 3 and 4, we didn't use the fact that these are sheaves.
\end{rmk}

\section{Fourier-Mukai transforms of varieties of $\kappa (X)=\dim X -1$}
 In this section we explain the important situation to which Theorem~\ref{MT}
 can be applied.
 Let us consider the situation when $K_X$ (or $-K_X$) is semi-ample, i.e.
 $|mK_X|$ is free for some $m >0$ (or $m<0$). When $K_X$ is semi-ample, 
we have the following morphism, called Iitaka fibration:
$$\pi _X \colon X\longrightarrow Z\cneq \Proj \bigoplus _{m\ge 0}
H^0 (X,mK_X).$$
Kodaira dimension of its generic fiber is zero. 
Let $Y\in FM(X)$ and 
$\Phi \colon D(X)\to D(Y)$ be an equivalence. 
Note that $K_Y$ is also semi-ample
by Corollary~\ref{support2}. By Corollary~\ref{cano}, 
the target of its Iitaka fibration is also $Z$. Let $\pi _Y \colon Y\to Z$ 
be the Iitaka fibration. Let us take a general closed point $p \in Z$. 
Let $X_p \cneq \pi _X ^{-1}(p)$, and $Y_p \cneq \pi _Y ^{-1}(p)$. 
Assume that kernel of $\Phi$ satisfies the condition 
as in Theorem~\ref{MT}, 
for example kernel of $\Phi$ is a sheaf. 
Then Theorem~\ref{MT} states that there exists equivalence
$\Phi _p \colon D(X_p)\to D(Y_p)$ such that the following diagram 
commutes:
$$(\diamondsuit) \qquad \begin{CD}
D(X) @> {\mathbf{L}i_p ^{\ast}}>>  D(X_p)  @> i_{p \ast}>>  D(X) \\
@V\Phi VV @V \Phi _p VV  @V \Phi VV \\
D(Y) @>{\mathbf{L}j_p^{\ast}}>>  D(Y_p) @>j_{p \ast}>> D(Y).
\end{CD}$$
Here $i_p$ and $j_p$ are inclusions, $i_p\colon 
X_p \hookrightarrow X$, $j_p\colon 
Y_p \hookrightarrow Y$. 
The conditions of kernels are satisfied if $\kappa (X)= \dim X-1$. 
Note that Fourier-Mukai partners of the varieties of $\kappa (X)=\dim X$ 
are studied in ~\cite{Ka1}.

\begin{thm}\label{mati}
Let $X$ be a smooth projective variety such that $K_X$ is semi-ample, 
and $\kappa (X)= \dim X -1$. 
Let $Y\in FM (X)$ and $\Phi \colon D(X)\to D(Y)$ be an equivalence. 
Then in the above notations, there exists 
an equivalence $\Phi _p \colon D(X_p)\to D(Y_p)$ such that the diagram 
$(\diamondsuit)$ commutes. 
\end{thm}
\begin{proof}
Let $\pP \in D(X\times Y)$ be a kernel of $\Phi$. It suffices to show 
$\pP \dotimes \oO _{X_p \times Y}$ is a sheaf, up to shift. 
Note that 
$$\pP \dotimes \oO _{X_p \times Y}\cong \pP \dotimes \oO _{X\times Y_p}, $$
by Step 2 of Theorem~\ref{MT}. By taking the functors whose kernels are left 
hand side, right hand side respectively, 
we can obtain the isomorphism of functors:
$$\Phi (\ast \dotimes \oO _{X_p})\cong \Phi (\ast)\dotimes \oO _{Y_p}.$$
Note that the above isomorphism can be also applied to derived categories 
of quasi-coherent sheaves. Let us consider $\Phi (\oO _x)$ for $x\in X_p$. 
Take a general morphism:
$$v_x \colon \Spec \mathbb{C}[[t_1, \cdots ,t_{d-1}]] \longrightarrow X, $$
which takes a closed point of $\Spec \mathbb{C}[[t_1, \cdots, t_{d-1}]]$
 to $x\in X_p$. Here $d\cneq \dim X$. 
Let $R_x \cneq v_{x\ast}\mathbb{C}[[t_1, \cdots, t_{d-1}]] \in \QCoh (X).$ 
Then $R_x \dotimes \oO _{X_p}
\cong \oO _x$, and 
\begin{align*}\Phi (\oO _x) &\cong \Phi (R_x)\dotimes \oO _{Y_p} \\
 &\cong j_{p\ast}\dL j_p ^{\ast}\Phi (R_x).
 \end{align*}
 Since $Y_p$ is one-dimensional, $\dL j_p ^{\ast}\Phi (R_x)$ is a direct 
 sum of its cohomologies. Since 
 $$\Hom _X (\oO _x ,\oO _x)\cong \Hom _Y (\Phi (\oO _x), \Phi (\oO _x))\cong 
 \mathbb{C}, $$
 we can conclude $\Phi (\oO _x)$ is a coherent $\oO _{Y_p}$-module, 
 up to shift. We may assume $\Phi (\oO _x)$ is a sheaf for general 
 $x \in X_p$. 
 Then for all $x\in X_p$, $\Phi (\oO _x)$ is a sheaf. Hence 
 $$\pP \dotimes \oO _{x\times Y} \cong \pP \dotimes \oO _{X_p \times Y}
 \dotimes p_1 ^{\ast}\oO _{R_x},$$
 is a sheaf. The above object is calculated by the spectral sequence:
 $$E_2 ^{p,q}=\tT or _{-p}^{\oO _{X\times Y}}(H^q (A), p_1 ^{\ast}\oO _{R_x})
 \Rightarrow H^{p+q}(\pP \dotimes \oO _{x\times Y}).$$
 Here $A\cneq \pP \dotimes \oO _{X_p \times Y}$. 
 The above spectral sequence degenerates at $E_2$-terms, since 
 $E_2 ^{p,q}=0$ for $p\le -2$.  Therefore if $k\neq 0$, 
 $H^k (A) \otimes p_1 ^{\ast}\oO _{R_x}=0$ for $x\in X_p$.
 This implies $H^k (A)=0$ for $k\neq 0$.  
 \end{proof}

As the immediate application, we generalize the theorem of Bondal and Orlov~\cite{B-O1}.
\begin{thm}
Let $C$ be an elliptic curve, $Z$ be a smooth projective variety.  
Assume that $K_Z$ or $-K_Z$ is ample. 
Then $FM(C\times Z)=\{ C\times Z \}$. 
\end{thm}
\begin{proof} 
We show the theorem when $K_Z$ is ample. The other case follows similarly. 
Let us take $Y\in FM(C\times Z)$, and let $\Phi \colon D(C\times Z)\to D(Y)$ be an equivalence.  
Since $C$ is an elliptic curve, the projection $C\times Z \to Z$ 
gives Iitaka fibration. 
Note that $K_Y$ is also semi-ample, and let $\pi \colon Y\to Z$ 
be its Iitaka fibration. 
Take a general closed point $p\in Z$ and fix it. 
Let $C^{\dag} \cneq \pi ^{-1}(p)$.
Then   
we can find an object $\uU \in D(C\times C^{\dag})$ such that 
 $\Phi _{C\to C^{\dag}}^{\uU} \colon D(C)\to D(C^{\dag})$ gives equivalence by Theorem~\ref{mati}. 
 Note that $C\cong C^{\dag}$, since Fourier-Mukai partners of a curve consists 
 of itself.  
On the other hand, as in Lemma~\ref{ob}, the following 
functor gives equivalence:
$${}\circ \uU  \colon D(C\times Z) \ni a \mapsto a \circ \uU  
\in D(C^{\dag} \times Z).$$
Let us compose the above equivalence with $\Psi \cneq \Phi ^{-1}$. We obtain the equivalence:
$$ ({}\circ \uU) \circ \Psi \colon D(Y) \longrightarrow D(C\times Z)
\longrightarrow D(C^{\dag}\times Z),$$
which takes $\oO _x$ to $\oO _{(x,p)}$ for all $x\in C^{\dag}$. 
Therefore we obtain the birational map over $Z$ by Lemma~\ref{bir} below,
$$f \colon Y \dashrightarrow C^{\dag} \times Z.$$
Note that $f$ is defined on the neighborhood of $C^{\dag}$. 
Since $Y$ and $C^{\dag} \times Z$ are both minimal models, 
$f$ is isomorphic in codimension one. 
We show that $f$ is in fact isomorphism. 
Let us take an ample divisor $H\subset Y$, and its strict transform 
$H^{\dag} \subset C^{\dag}\times Z$. 
It suffices to show $H^{\dag}$ is nef. But this is clear since
$H^{\dag}$ is effective, and we can deform $H^{\dag}$ freely using translations
of $C^{\dag}$. \end{proof}

\section{Fourier-Mukai partners of 3-folds of $\kappa (X)=2$}
In this section, we study $FM(X)$ when $\dim X=3$ and $\kappa (X)=2$. 
The relative moduli spaces of stable sheaves for three-dimensional Calabi-Yau fibrations 
are studied in~\cite{B-M2}. Combining Theorem~\ref{MT} with 
their results, we can study $FM(X)$ in this case. 
Before that, we recall some terminology of birational geometry, 
and give some useful lemmas.

\begin{defi}
Let $X$ and $Y$ be projective varieties with only canonical singularities.
 A birational map 
$\alpha \colon 
X \dashrightarrow Y$ is called crepant, if there exists a smooth projective
variety $Z$ and birational morphisms $f\colon Z \to X$, $g\colon Z\to Y$, 
such that $\alpha \circ f=g$, and $f^{\ast}K_X =g^{\ast}K_Y$. 
In this case, we say $X$ and $Y$ are $K$-equivalent under $\alpha$. 
\end{defi}

The following birational transform called ``flop" is a special kind of crepant birational map. 

\begin{defi}
Let $X$ and $Y$ be projective varieties with only canonical singularities. 
A birational map $\alpha \colon X \dashrightarrow Y$ is called a flop, 
if there exist a normal projective variety $W$ and crepant 
birational morphisms $\phi \colon X\to W$, $\psi \colon Y\to W$ which 
satisfy the following:
\begin{itemize}
\item $\phi =\psi \circ \alpha$.
\item $\phi$ and $\psi$ are isomorphisms in codimension one. 
\item Relative Picard numbers of $\phi$, $\psi$ are one. 
\item Let $H$ be a $\phi$-ample divisor on $X$, and $H'$ be 
its strict transform on $Y$. Then $-H'$ is $\psi$-ample. 
\end{itemize}
\end{defi}

In this paper, we will only use flops of smooth 3-folds. In dimension three, 
crepant birational maps are connected by finite number of 
flops~\cite[Theorem 4.6]{Ka1}. 
Next we give some useful lemmas.

\begin{lem}\label{bir}
Let $X$ and $Y$ be smooth projective varieties, and 
$\Phi \colon D(X)\to D(Y)$ be an equivalence. Assume 
for some closed point $x\in X$, we have $\dim \Supp \Phi (\oO _x)=0$.
Then there exists an open neighborhood $U$ of $x$, and
$r\in \mathbb{Z}$,
such that for $x' \in U$, there exists $f(x') \in Y$ which satisfies 
$\Phi (\oO _{x'})=\oO _{f(x')}[r]$. Moreover $X$ and $Y$ are $K$-equivalent 
under birational map $f \colon X \dashrightarrow Y$. \end{lem}
\begin{proof}
Since $\Phi$ gives an equivalence, we have
$$\Ext _Y ^i (\Phi (\oO _x), \Phi (\oO _x))= \left\{
\begin{array}{cc} 0 & (i<0), \\ \mathbb{C} & (i=0). 
\end{array}\right.$$
Then using the same argument as in~\cite[Proposition 2.2]{B-O1},
 there exists a point $y\in Y$
and $r\in \mathbb{Z}$ such that $\Phi (\oO _x) \cong \oO _y [r]$. 
Then as in~\cite[Theorem 2.5]{B-M1}, we can find a desired $U$ and a 
birational map $f$. Let $\pP \in D(X\times Y)$ be a kernel of $\Phi$. 
Since 
$$\pP \otimes p_1 ^{\ast}\omega _X \cong \pP \otimes p_2 ^{\ast}\omega _Y, $$
as in Section 4, $p_1 ^{\ast}\omega _X$ and $p_2 ^{\ast}\omega _Y$ 
are numerically equal on $\Supp \pP$. By the construction of $f$, 
we have $\Gamma _f \subset \Supp \pP$, where $\Gamma _f$ is a graph 
of $f$. Therefore $X$ and $Y$ are $K$-equivalent under $f$. 
\end{proof}

\begin{lem}\label{Bs}
Let $X$ and $Y$ be smooth projective varieties, and $\Phi \colon 
D(X) \to D(Y)$ gives equivalence. Then 
the followings hold:

(i) For a closed point $x\in X$, $\omega _Y$ is numerically zero on 
$\Supp \Phi (\oO _x)$. 

(ii) If $x\in \Bs |mK_X|$, then $\Supp \Phi (\oO _x)\subset \Bs |mK_Y|$. 

(iii) If $x\notin \Bs |mK_X|$, then $\Supp \Phi (\oO _x) \cap \Bs |mK_Y|=\emptyset$.
\end{lem}
\begin{proof}
(i) Since $\Phi$ and Serre functor commutes, we have 
$\Phi (\oO _x)\otimes \omega _Y \cong \Phi (\oO _x)$.
(i) follows from this. 

(ii) This follows from Lemma \ref{support} immediately. 

(iii) Take $x\notin \Bs |mK_X|$ and assume that there exists $y\in \Supp \Phi (\oO _x)\cap \Bs |mK_Y|$. Then there exists a non-zero map 
$\Phi (\oO _x) \to \oO _y [i]$ for some $i$. Therefore there exists a non-zero map $\oO _x \to \Psi (\oO _y)[i]$. Since $\Psi (\oO _y)[i]$ is supported on 
$\Bs |mK_X|$, this is a contradiction. \end{proof}

 Now we state our main theorem of this section. 
\begin{thm}\label{k1}Let $X$ be a smooth projective 3-fold of $\kappa (X)=2$. Then $Y\in FM(X)$ if and only if 
one of the following holds:

(1) $X$ and $Y$ are connected by finite number of flops. 

(2) There exists a following diagram:

$$\xymatrix{
Y \ar@{.>}[r]^{\emph{flops}}  & J^H (d) \ar[rd] _{\pi '}&  & M \ar[ld]^{\pi} 
\ar@{<.}[r]^{\emph{flops}}&X\\
               &                & S, & &\\
}$$
where $\pi \colon M\to S$ is an elliptic fibration with $\omega _{M} \equiv _{\pi}0$, $H\in \Pic (M)$ is a polarization, and
$d\in \mathbb{Z}$. $J^{H}(d)\subset M^{H}(M/S)$ is an irreducible component which is fine, and contains a point corresponding to line bundles of degree $d$ on 
smooth fibers of $\pi$.  

\end{thm}
 ``If" direction is already proved in ~\cite[Theorem 8.3]{B-M2}
  and ~\cite{Br1}, when $S$ is smooth. 
  We can check that the assumption ``$S$ is smooth" is not required in 
  their proof, hence the ``if" direction holds. 
  We prove ``only if" direction. Let us take $Y\in FM(X)$.
We use the same notations as in the previous sections. In particular
$\Phi \colon D(X)\to D(Y)$ gives equivalence, $\pP \in D(X\times Y)$ 
is a kernel of $\Phi$, and $\Psi$ is a quasi-inverse of $\Phi$. 
Note that, by Lemma~\ref{bir}, we may assume $\dim \Supp \Phi (\oO _x)\ge 1$ 
for all $x\in X$. In this situation, we will construct a diagram (2), or 
show (1) holds. We divide the proof into 5 Steps. 

\begin{step}
Application of Theorem~\ref{MT}. \end{step}
At first, we apply Theorem~\ref{MT}, and give the preparation for 
the proof. Since $\dim X=\dim Y=3$, we can run minimal 
model programs, and obtain birational minimal models $X_{\textrm{min}}$ and $Y_{\textrm{min}}$:
$$\xymatrix{
X \ar@{.>}[r]^{\textrm{MMP}}  & X_{\textrm{min}} \ar[rd] _{\pi _X}&  
& Y_{\textrm{min}} \ar[ld]^{\pi _Y} 
\ar@{<.}[r]^{\textrm{MMP}}&Y\\
               &                & Z. & &\\
}$$
Here $\pi _X$, $\pi _Y$ are Iitaka fibrations.  Note that $\dim Z =2$, and 
generic fibers of $\pi _X$, $\pi _Y$ are elliptic curves.
Then for sufficiently large $m$, we obtain isomorphisms,
$$X\setminus \Bs |mK_X|\stackrel{\cong}{\longrightarrow}X_{\textrm{min}}\setminus C_X, \qquad
Y\setminus \Bs |mK_Y|\stackrel{\cong}{\longrightarrow}Y_{\textrm{min}}\setminus C_Y,$$
for some closed subsets $C_X \subset X_{\textrm{min}}$, $C_Y \subset Y_{\textrm{min}}$ with $\dim C_X \le 1$, $\dim C_Y \le 1$. 
Let us take general members $E_i \in |mK_X|$, for $i=1,2$. 
By Corollary~\ref{cano}, we have the isomorphism of linear systems:
$$|mK_X| \stackrel{\sim}{\longrightarrow} |mK_Y|.$$
Let $E_i ^{\dag}\in |mK_Y|$ corresponds to $E_i$.  Also note that 
we have the isomorphisms:
$$H^0 (X,mK_X) \cong H^0 (X_{\textrm{min}}, mK_{X_{\textrm{min}}}), \quad 
H^0 (X,mK_Y) \cong H^0 (Y_{\textrm{min}}, mK_{Y_{\textrm{min}}}), $$
for sufficiently divisible $m$. Let 
$$E_i '\in |mK_{X_{\textrm{min}}}|, \quad E_i ^{' \dag}\in 
|mK_{Y_{\textrm{min}}}|,$$
correspond to $E_i$, $E_i ^{\dag}$ under the above isomorphisms respectively. 
Then if we choose $E_i$ sufficiently general, then we have,
$$E_1 ' \cap E_2 ' \cap C_X =\emptyset.$$ 
Therefore we have the following decompositions:
$$E_1 \cap E_2  =(E_1 ' \cap E_2 ') \coprod \Bs |mK_X|, \qquad E_1 ^{\dag} \cap E_2 ^{\dag} =(E_1^{' \dag} \cap E_2 ^{' \dag}) \coprod \Bs |mK_Y|.$$
Now let us take $C\in \pi _0 (E_1 ' \cap E_2 ')$. 
We can consider $C$ as 
a curve on $X$. Using Corollary~\ref{support2} and Lemma~\ref{Bs}, we can find 
$C^{\dag}\in \pi _0 (E_1 ^{'\dag}\cap E_2 ^{'\dag})$ such that $\Phi$ takes 
$D_C (X)$ to $D_{C^{\dag}}(Y)$. 
Now using the same argument as in Theorem~\ref{mati}, we can see $\pP \dotimes 
\oO _{C\times Y}$ is a sheaf, up to shift. Then we can apply Theorem~\ref{MT}, 
so there exists an equivalence $\Phi _{C}\colon D(C)\to D(C^{\dag})$, 
such that the diagram of Theorem~\ref{MT} commutes. 

\begin{step} Construction of $M$. \end{step}
In this step, we will construct a desired $M$. We will 
construct $M$ as a moduli space of stable sheaves on $Y$.
Let us take $x\in C$, and consider $\Phi _C (\oO _x)\in D(C^{\dag})$.
As in Theorem~\ref{mati}, we may assume 
$\Phi _C (\oO _x)$ is a simple sheaf on $C^{\dag}$. 
Since $C^{\dag}$ is an elliptic curve, $\Phi _C(\oO _x)$ is a stable 
sheaf on $C^{\dag}$. 
 Let $\rk \Phi _C (\oO _x)=a$ and $\deg \Phi _C (\oO _x)=b$. 
By the commutative diagram of Theorem~\ref{MT}, $\Phi (\oO _x)$ is a 
stable sheaf on $Y$ supported on $C^{\dag}$, with respect to any polarization. 
Then take a polarization $H' \in \Pic (Y)$, and 
consider moduli space of stable sheaves $M^{H'}(Y/\Spec \mathbb{C})$.
Let 
$$M \subset M^{H'}(Y/\Spec \mathbb{C})$$
 be an irreducible component, which contains a point corresponding to 
 $\Phi (\oO _x) \in \Coh (Y)$. 
 Note that there exists a birational map:
 $$f_1 \colon X \dashrightarrow M,$$
 which takes a general point $x\in X$ to a point of $M$, corresponding to a
 stable sheaf $\Phi (\oO _x)$.
 We show $M$ is a fine moduli scheme, or (1) holds. 
  For $E, F\in D(X)$, we define $\chi (E, F)$ as follows:
$$\chi (E, F)\cneq \sum (-1)^i \dim \Ext _X ^i (E,F).$$
Since $\chi (\Phi (\oO _X), \Phi (\oO _x))=\chi (\oO _X, \oO _x)=1$, Riemann-Roch implies that 
$$b\cdot \ch _0 \Phi (\oO _X)^{\ast}+a (c_1 (\Phi (\oO _X)^{\ast})\cdot C^{\dag})=1.$$
Here $\Phi (\oO _X)^{\ast}$ means derived dual of $\Phi (\oO _X)$. 
We divide into 2-cases:
\begin{ccase} $b=0$ \end{ccase}
If $b=0$, then $a=c_1 (\Phi (\oO _X)^{\ast})\cdot C^{\dag}=1$. 
Therefore there exists an effective divisor $E$ on $Y$ such that 
$E\cdot C^{\dag}=1$.  There exists a birational map:
$$f _2 \colon Y\dashrightarrow M, $$
such that $f _2 $ takes general 
point $y\in Y$ to a point corresponding to $\oO _{C_y}(E\cap C_y -y)$, 
a degree zero line bundle on $C_y$. Here  
$C_y$ is a compact fiber of the Iitaka fibration $Y\dashrightarrow Z$, which contains $y$. Composing these we obtain a birational map, 
$$f\cneq f _2 ^{-1}\circ f _1 \colon X\dashrightarrow M\dashrightarrow 
Y,$$ which satisfies
 $f (x)\in \Supp \Phi (\oO _x)$ for general $x\in X$. 
Therefore $\Gamma _{f}\subset \Supp \pP$, where $\Gamma _{f}$ is a graph of 
$f$. Since $p_1^{\ast}K_X \equiv p_2 ^{\ast}K_Y$ on 
$\Supp \pP$, it is also true on $\Gamma _{f}$. Therefore $X$ and $Y$ are $K$-equivalent under birational map $f$. 

\begin{ccase}$b\neq 0$ \end{ccase}
Let us replace $H'$ to $\det \Phi (\oO _X)^{\ast}\pm lbH'$ for $l\gg 0$. Then
 we may assume $\ggcd (a(H' \cdot C^{\dag}), b)=1$. 
Then 
\begin{align*}
\ggcd \{ \chi (\Phi (\oO _x)\otimes H ^{'\otimes m})\mid m\in \mathbb{Z}\}
 &= \ggcd
  \{ ma (H' \cdot C^{\dag})+b \mid m\in \mathbb{Z} \} \\
&= 1. 
\end{align*}
By Lemma~\ref{fine}, this implies that $M$ is a fine moduli scheme. 

\begin{step}
$M$ is smooth and the universal sheaf $\uU \in \Coh (Y\times M)$ gives an equivalence 
$$\Phi _M \cneq \Phi _{M\to Y}^{\uU}\colon D(M)\longrightarrow
 D(Y).$$ \end{step}
\begin{proof}
For $p \in M$, let $\uU _p \in \Coh (Y)$ be the corresponding stable sheaf. 
We check that the conditions of Theorem~\ref{int} are satisfied. 
First we show $\uU _p \otimes \omega _Y \cong \uU _p$. 
Let 
$$\xymatrix{
 \widetilde{X} \ar[d]_{g} \ar[rd] ^{h}&  \\ 
X  \ar@{.>}[r]^{f_1}  & M \\
}$$
be an elimination of indeterminacy. Consider morphisms,
$$g\times \id \colon \widetilde{X}\times Y \longrightarrow X\times Y, \quad
h\times \id \colon \widetilde{X}\times Y \longrightarrow M\times Y, $$
and objects, 
$$\dL (g \times \id )^{\ast}\pP \in D(\tX \times Y), \quad (h \times 
\id)^{\ast}\uU \in \Coh (\tX \times Y).$$ 
Take $x\in \tX$ and let $i_{{x}\times Y}\colon x \times Y \hookrightarrow X\times Y$ be an inclusion. Then 
\begin{align*}\dL i_{x \times Y}^{\ast}\circ \dL (g \times \id)^{\ast}\pP &=\dL i_{g(x) \times Y}^{\ast}\pP =\Phi (\oO _{g (x)}) \\
\dL i_{x \times Y}^{\ast}\circ (h \times \id)^{\ast}\uU &= \uU _{h(x)}. 
\end{align*}
Take open subsets $X^0 \subset X$, $Y^0 \subset Y$, $Z^0 \subset Z$ such that the rational maps $X\dashrightarrow Z$, $Y\dashrightarrow Z$ are defined on $X^0$, 
$Y^0$, and $X^0 \to Z^0$, $Y^0 \to Z^0$ are smooth projective.
From here, we will shrink $Z^0$ if necessary. 
 Since $f_1$ is defined on $X^0$, we can think $X^0$ as an open subset of $\widetilde{X}$. So if $x\in X^0 \subset \tX$, then $\Phi (\oO _{g (x)})=\uU _{h (x)}$.
 This implies 
$$\Supp (h \times \id)^{\ast}\uU \cap (X^0 \times Y)=\Supp \dL (g \times \id)^{\ast}\pP \cap (X^0 \times Y).$$
Therefore by Lemma~\ref{irr} below, we have 
$$\Supp (h \times \id)^{\ast}\uU \subset 
\Supp \dL (g \times \id)^{\ast}\pP \subset \widetilde{X} \times Y.$$
Therefore for all ${x}\in \tX$, we have 
$$\Supp (h \times \id)^{\ast}\uU \cap (x \times Y) \subset \Supp \dL (g \times \id)^{\ast}\pP \cap (x \times Y).$$
So $\Supp \uU _{h ({x})}\subset \Supp \Phi (\oO _{g ({x})})$
 follows.
Since $\omega _Y$ is numerically zero on $\Supp \Phi (\oO _{g ({x})})$, this is also true on $\Supp \uU _{h ({x})}$, hence on 
$\Supp \uU _p$ for all $p\in M$. Therefore $\uU _p \otimes \omega _Y$ is also $H'$-stable, and its reduced Hilbert polynomial is equal to $\uU _p$, i.e. 
$$p(\uU _p , H')=p(\uU _p \otimes \omega _Y , H').$$ 
On the other hand, there exists a non-trivial map $\uU _p \to \uU _p \otimes \omega _Y$ by semi-continuity. So $\uU _p \cong \uU _p \otimes \omega _Y$ for all $p\in M$. 

Secondly we show the set 
$$\Gamma (\uU )\cneq \{ (p_1, p_2)\in M\times M \mid \Ext _Y ^i (\uU _{p_1}, \uU _{p_2})\neq 0 \quad \textrm{for some }i\in \mathbb{Z}\}$$
has $\dim \Gamma (\uU)\le 4$. It suffices to show if $(p_1 , p_2)\in \Gamma (\uU)\setminus \Delta _M$, where $\Delta _M$ is a diagonal, then 
$p_i \in M\setminus f_1 (X^0)$. Assume that $p_1 \in f_1 (X^0)$. Since 
$\Ext _Y ^i (\uU _{p_1}, \uU _{p_2})\neq 0$, we have $\Supp \uU _{p_1}\cap \Supp \uU _{p_2}\neq \emptyset$. Take an irreducible component 
$l\subset \Supp \uU _{p_2}$ such that $\Supp \uU _{p_1} \cap l \neq \emptyset$. Since we have assumed $p_1 \in f_1 (X^0)$, we have 
$$\Supp \uU _{p_1}\cap \Bs |mK_Y|=\emptyset.$$
So it 
follows that $l$ is not contained in $\Bs |mK_Y|$. Furthermore 
$K_Y \cdot l =0$, since $K_Y$ is numerically zero on $\Supp \uU _p$. Therefore
 $l \cap \Bs |mK_Y|=\emptyset$ and 
$l$ is contained in the fiber of the Iitaka fibration, 
$Y\setminus \Bs |mK_Y|\to Z$. This implies $l=\Supp \uU _{p_1}$ and therefore $\Supp \uU _{p_2}=\Supp \uU _{p_1}$, since
$\Supp \uU _{p_2}$ is connected. Therefore $\uU _{p_2}$ is a stable sheaf on $\Supp \uU _{p_1}$, so $p_2 \in f_1 (X^0)$. Let $q_i \in X^0$ be 
points such that $p_i =f_1 (q_i)$. Then 
$\Ext _Y ^i (\uU _{p_1}, \uU _{p_2})=\Ext _X ^i (\oO _{q_1}, \oO _{q_2})\neq 0$ implies $q_1 =q_2$ and $p_1 =p_2$. But this contradicts to $(p_1, p_2) \notin \Delta _M$. 
\end{proof} 
In the above proof, we used the following lemma:
\begin{lem}\label{irr}
$\Supp (h \times \emph{id} )^{\ast}\uU $ is irreducible.\end{lem} 
\begin{proof}
Let $\widetilde{f}\colon \tX \times Y \to \tX$ be a projection.
Note that a general fiber of the restriction of $\widetilde{f}$ to $\Supp (h \times \id)^{\ast}\uU$ is an elliptic curve. Therefore if $\Supp (h \times \id)^{\ast}\uU$ is not irreducible, 
then there exists $p\in \Ass ((h \times \id)^{\ast}\uU)$ such that $\dim \oO _{\tX ,\widetilde{f}(p)}\ge 1$. Take a non-zero element of the maximal ideal 
$t\in m_{\widetilde{f}(p)}\subset \oO _{\tX ,\widetilde{f}(p)}$. Then $\oO _{\tX ,\widetilde{f}(p)}\stackrel{\times t}{\to}\oO _{\tX ,\widetilde{f}(p)}$ 
is injective. Since $(h \times \id)^{\ast}\uU$ is flat over $\tX$, we have an injection,
$$((h \times \id)^{\ast}\uU)_p \stackrel{\times \widetilde{f}^{\ast}t}{\longrightarrow}((h \times \id)^{\ast}\uU)_p,$$
and $\widetilde{f} ^{\ast}t\in m_p \oO_{\tX \times Y, p}$. But this contradicts to $p\in \Ass ((h \times \id)^{\ast}\uU)$. \end{proof}
\begin{step}
$X$ and $M$ are connected by finite number of flops, and $M$ has an 
elliptic fibration $\pi \colon M \to S$ with $\omega _M \equiv _{\pi}0$.
\end{step}

\begin{proof}
Consider the following composition:
$$\Psi \circ \Phi _{M}\colon 
D(M)\longrightarrow D(Y) \longrightarrow D(X).$$
This is an equivalence and for general points $p\in M$, 
we have 
$$\dim \Supp \Psi \circ \Phi _M(\oO _p)=0.$$
Therefore $X$ and $M$ are connected by finite number of flops. 
Since $\Supp \uU \subset Y\times M$ is irreducible 
and all the fibers of the projection 
$\Supp \uU \to M$ are one-dimensional, this is 
a well-defined family of proper algebraic cycles in the sense of 
~\cite{Ko}. Therefore there exists a morphism $M\to \Chow(Y)$ which 
takes $p\in M$ to an algebraic cycle whose support is equal to $\Supp \uU _p$.  
Let 
$$M\stackrel{\pi}{\to} S \to \Chow (Y)$$
be a stein factorization. We show that $\omega _{M}\equiv _{\pi}0$. 
Let us take $p,p'\in M$ such that $\pi (p)=\pi '(p)$. Then by the definition of $\pi$, 
it follows that 
$$\Supp \Phi _M(\oO _p)=\Supp \Phi _{M}(\oO _{p'}).$$
Let us take $q\in \Supp \Phi _{M}(\oO _p)$. Then 
$p' \in \Supp (\Phi _M)^{-1}(\oO _q)$. Therefore $\pi ^{-1}\pi (p)\subset \Supp (\Phi _{M})^{-1}(\oO _q)$. 
This implies $\omega _{M}\equiv _{\pi} 0$. 
\end{proof}

\begin{step}There exists a polarization $H\subset M$, $d\in \mathbb{Z}$, 
such that $J^H (d)$ is fine, and smooth. Moreover
 $Y$ and $J^H (d)$ are connected 
by finite number of flops. \end{step}
\begin{proof}
We continue the same argument. 
Let us take a general closed point $y\in Y$. The object
$$(\Phi _{M})^{-1}(\oO _y) \in D(M)$$
is a stable sheaf on a general fiber of $\pi$. Let its rank and degree be $c$ and $d$ respectively. 
Let $H\in \Pic (M)$ be a polarization, and take an irreducible component $M^{+}\subset M^{H}(M/S)$
 which contains a point corresponding to $(\Phi _M)^{-1}(\oO _y)$.
Similarly, take an irreducible component $J^{H}(d)\subset M^{H}(M/S)$ which contains a point corresponding to 
line bundles of degree $d$ on smooth fibers of $\pi$.  
By the same argument as before, we can choose $H$ such that 
$$\pi{''}\colon M^{+}\to S, \quad \pi '\colon J^{H}(d)\to S, $$
are fine moduli schemes (or $X$ and $Y$ are connected by finite 
number of flops if $d=0$). By ~\cite[Theorem 8.3]{B-M2}, $M^{+}$ and $J^{H}(d)$ are smooth, $\omega _{M^{+}}\equiv _{\pi ''}0$, 
$\omega _{J^{H}(d)}\equiv _{\pi '}0$, and the universal sheaf $\vV \in \Coh (M^{+}\times _S M )$ gives an equivalence
$$\Phi _{M^{+}}\cneq \Phi _{M^{+}\to M}^{\vV} \colon D(M^{+})\stackrel{\sim}{\longrightarrow}
 D(M).$$
Since the composition 
$$\Phi _{M}\circ \Phi _{M^{+}}\colon D(M^{+})\stackrel{\sim}{\longrightarrow} D(M)\stackrel{\sim}{ 
\longrightarrow} D(Y)$$
takes general points to general points, $Y$ and $M^{+}$ are connected by finite number of flops. By ~\cite[Theorem 6]{A}, 
there exists a following birational map over 
$S$:
$$M^{+} \ni E \mapsto \wedge ^c E \in J^{H} (d). $$
Since they are both minimal over $S$, $M^{\dag}$ and $J^{H} (d)$ are connected by finite number of flops. Now we have obtained the diagram (2). 
\end{proof}

If $X$ is minimal we have a better result. By the abundance theorem in dimension three, $K_X$ is semi-ample. 
Let $\pi _X \colon X\to Z$ be its Iitaka fibration. We define $\lambda _X >0$ as follows:
$$\lambda _X \cneq \ggcd \{ c_1 (E)\ccdot f_X \mid E \in D(X) \}, $$
where $f_X$ is a cohomology class of a general fiber of $\pi _X$. 
For a polarization $H$ on $X$, let  
$J^H (b) \subset M^H (X/Z)$
be as in the Theorem \ref{k1}.
The proof of the following theorem is almost the same as in the previous theorem and is
left to the reader.

\begin{thm}
Let $X$ be a smooth minimal 3-fold with $\kappa (X)=2$. Then 
$Y \in FM(X)$ if and only if there exists some $b\in \mathbb{Z}$ which is co-prime to $\lambda _X$, 
and there exists a polarization $H$ on $X$, for which $J^H (b)$ is a fine moduli scheme, such 
that Y and $J^H (b)$ are connected by finite number of flops.
\end{thm}

Because $J^H (b+\lambda _X)\cong J^H (b)$, birational classes of $FM(X)$ are finite in the above case. By ~\cite{Ka2}, the number of 3-dimensional minimal model in 
a fixed birational class is finite. So we obtain the next corollary.  
\begin{cor}
Let $X$ be a smooth minimal 3-fold with $\kappa (X)=2$. Then $\sharp FM(X)< \infty$. 
\end{cor}

{\bf acknowledgements} 
The author would like to express his profound gratitude to 
Professor Yujiro Kawamata, for many valuable comments, and warm encouragement.

\end{document}